\documentclass[11pt]{article}
\usepackage{inputenc}
\usepackage[T2A]{fontenc}
\usepackage[dvips]{graphicx}
\usepackage{amsthm,amsfonts,amsmath,amssymb}
\usepackage{color}
\usepackage{amsmath,bm}

\begin{document}

\author{\bf A.V. Rukavishnikov }
\title{\bf New approach for solving stationary  nonlinear Navier-Stokes equations in non-convex domain }
\maketitle{}

	\centerline{Computing Center of the Far Eastern Branch of the Russian Academy of Sciences,}
\centerline{680000, Khabarovsk, Russia, e-mail: \tt 78321a@mail.ru}			
\maketitle{}
\begin{abstract} \begin{tabular}{p{0mm}p{139mm}}
&\noindent {\footnotesize \qquad In the paper, an approach for the numerical solution of stationary
 nonlinear Navier-Stokes equations in rotation and convective forms in a polygonal domain
 containing one reentrant corner on its boundary, that is, a corner greater than $ \pi $ is considered.
The method allows us to obtain the 1st order of convergence of the approximate solution to the exact one with respect to the grid step h,
 regardless of the reentrant corner value.

{\bf Keywords:} {nonlinear Navier-Stokes equations, singularity, finite element method.}}
\end{tabular}\end{abstract}

\paragraph{Introduction.}

In the paper, an approach for the numerical solution of stationary
 nonlinear Navier-Stokes equations in rotation and convective forms in a polygonal domain
 containing one reentrant corner on its boundary, that is, a corner greater than $ \pi $ is considered.
 The problem is linearized using the Picard iterative procedure.
Taking into account the principle of coordinated estimates,
the approximate solution of such problem,
obtained using classical finite-difference and finite-element schemes,
converges to the exact one with a rate significantly less than 1st with respect to the grid step h
and depends on the value of the reentrant corner, i.e. the  pollution effect takes place \cite{Blum}.
The approximate approach proposed by us is based on the definition of the $ R_{\nu}$-generalized
solution and the introduction of auxiliary weight functions in the basis of finite-dimensional spaces.
The method allows us to obtain the 1st order of convergence of the approximate solution to the exact one with respect to the grid step h,
 regardless of the reentrant corner value.

For the first time (see \cite{Ruk1989}) it was proposed to define the solution as $ R_{\nu} $-generalized one
 for elliptic boundary value problems. It has been proved (see, for example, \cite{Ruk2014})
 its existence and uniqueness in weighted spaces and sets.
It was established that the $ R_{\nu} $-generalized solution (under certain conditions)
belongs to weighted spaces of a higher order \cite{Ruk2009}.
Numerical methods for solving elliptic boundary value problems
have been developed and  the rate of convergence have been obtained (see \cite{Ruk2000}).
The weighted finite element method for problems of elasticity theory with a corner singularity and a crack was proposed in \cite{Ruk2015} and \cite{Ruk2021}, respectively,
for the Stokes and Oseen problem in various forms in \cite{Ruk2018} and \cite{Ruk2019}, \cite{Ruk2020}, respectively.

In the presented paper, a series of numerical experiments was carried out for various values
of the free parameters of the method for both forms. Some of the free parameters
 do not depend on the value of the reentrant corner for both forms of the nonlinear problem.
The advantage of the proposed method over classical approaches has been established.
Finally, conclusions are drawn.

\paragraph{1. The problem statement. Definition of the $R_{\nu}$-generalized solution.}

We will consider a nonlinear problem obtained by applying an implicit time scheme to incompressible
 Navier-Stokes equations
 for a viscous fluid flow in convective form in two forms a two-dimensional domain $\Omega$ of Euclidean space
(${\bf x}=(x_1, x_2) $ is its element) with the boundary
$\partial \Omega, \bar{\Omega}=\Omega\cup \partial \Omega.$

{\bf Problem 1} (the convective form of Navier-Stokes equations). Find the velocity vector function
${\bf w}={\bf w}({\bf x})=(w_1({\bf x}), w_2({\bf x}))$ and kinematic pressure
 $\bar q=\bar q({\bf x})$, satisfying the system of equations
\begin{gather}
\alpha\, {\bf w}-{\mu}\triangle {\bf w} +({\bf w} \cdot \nabla) {\bf w}+ \nabla \bar q={\bf f}\qquad \,\,\,\,\,\mbox{ in } \,\,\qquad \Omega,
\label{eq:1}\\
\,\,\,\qquad\qquad\qquad\qquad\qquad\,\,\mbox{ div } {\bf w} ={ 0} \qquad \mbox{ in } \qquad \,\,\Omega
\label{eq:2}
\end{gather}
and boundary conditions
\begin{gather}
\,\,\,\qquad\qquad\qquad\qquad\qquad\qquad\,\,\,\,\, {\bf w} =\bold{\sigma} \qquad \mbox{ on } \qquad \partial \Omega,
\label{eq:3}
\end{gather}
where $ {\bf f} = {\bf f} ({\bf x}) $ and $ \bold{\sigma} = \bold{\sigma} ({\bf x}) $ are known functions
 of the right-hand sides in $ \Omega $ and on
$ \partial \Omega $, respectively, $ \alpha $ is a positive discretization parameter,
$ \mu = \frac {1}{Re}, Re$ is the Reynolds number.

{\bf Problem 2} (the rotation form of Navier-Stokes equations). Find the velocity vector function
${\bf w}={\bf w}({\bf x})=(w_1({\bf x}), w_2({\bf x}))$ and Bernoulli pressure
 $\tilde q=\tilde{q}({\bf x})$, satisfying the system of equations and boundary conditions
\begin{gather}
\alpha\, {\bf w}-{\mu}\triangle {\bf w} +\mbox{ curl }{\bf w} \times {\bf w}+ \nabla \tilde{q}={\bf f}\qquad \mbox{ in } \,\,\qquad \Omega,
\label{eq:4}\\
\,\,\,\qquad\qquad\qquad\qquad\qquad\,\,\mbox{ div } {\bf w} ={ 0} \qquad \mbox{ in } \qquad \,\,\Omega,
\label{eq:5}\\
\,\,\,\qquad\qquad\qquad\qquad\qquad\qquad\,\,\,\,\, {\bf w} =\bold{\sigma} \qquad \mbox{ on } \qquad \partial \Omega.
\label{eq:6}
\end{gather}

In the equations (\ref{eq:1})-(\ref{eq:6}) and further, we use the notation:
$\triangle {\bf z}=(\triangle z_1, \triangle z_2)^T, \triangle z_i=\frac{\partial^2 z_i}{\partial x_1^2}+
\frac{\partial^2 z_i}{\partial x_2^2}, \mbox{ curl }{\bf z}=\frac{\partial z_2}{\partial x_1}-
\frac{\partial z_1}{\partial x_2} \mbox{ div }{\bf z}=\frac{\partial z_1}{\partial x_1}+
\frac{\partial z_2}{\partial x_2}, s\times {\bf z}=\Bigl(-s z_2, s z_1\Bigr)^T
({\bf z}\cdot \nabla) {\bf v}=\Bigl(z_1 \frac{\partial v_1}{\partial x_1}+
z_2 \frac{\partial v_1}{\partial x_2}, z_1 \frac{\partial v_2}{\partial x_1}+
z_2 \frac{\partial v_2}{\partial x_2}  \Bigr)^T$
and $\nabla s=\Bigl(\frac{\partial s}{\partial x_1}, \frac{\partial s}{\partial x_2}\Bigr)^T$ respectively.

{\bf Remark 1}.
{\it Problems (\ref{eq:1})-(\ref{eq:3}) and (\ref{eq:4})-(\ref{eq:6}) are obtained from each other using identity}
$
 ({\bf w} \cdot \nabla) {\bf w} =(\mbox{ curl }{\bf w})\times {\bf w}+\frac{1}{2}\nabla {\bf w}^2
$
{\it and equality}
$
\tilde{q}=\bar{q}+ \frac{1}{2} {\bf w}^2.
$

Next, we write Problem 1 ana Problem 2 together, assuming that for $\gamma=1$ we obtain a convective form of Navier-Stokes equations and for $\gamma=0$ --- a rotation form of nonlinear problem.

{\bf Problem 3} (joint formulation of a nonlinear problem). Find the velocity vector function
${\bf w}={\bf w}({\bf x})=(w_1({\bf x}), w_2({\bf x}))$ and pressure
 $q=q({\bf x})$, satisfying the system of equations and boundary conditions
\begin{gather}
\alpha\, {\bf w}-{\mu}\triangle {\bf w} +\gamma({\bf w} \cdot \nabla) {\bf w}+(1-\gamma)\mbox{ curl }{\bf w} \times {\bf w}+ \nabla \tilde{q}={\bf f}\qquad \mbox{ in } \,\,\qquad \Omega,
\label{eq:7}\\
\,\,\,\qquad\qquad\qquad\qquad\qquad\qquad\qquad\qquad\qquad\quad\,\,\,\,\,\mbox{ div } {\bf w} ={ 0} \qquad \mbox{ in } \qquad \,\,\Omega,
\label{eq:8}\\
\,\,\,\qquad\qquad\qquad\qquad\qquad\qquad\qquad\qquad\qquad\quad\qquad\,\,\,\,\,\,\, {\bf w} =\bold{\sigma} \qquad \mbox{ on } \qquad \partial \Omega.
\label{eq:9}
\end{gather}

Since the system of {\bf Problem 3} is nonlinear, it should be linearized.
It is well known that if the Reynolds number and the function $ {\bf f} $ are not so large, then
the nonlinear problem (\ref{eq:7})-(\ref{eq:9}) has a unique solution $ ({\bf w},q) $.
   Let's construct an iterative Picard procedure
such that the sequence of iterative solutions $ ({\bf w}_k, q_k) $ of linear systems
\begin{gather}
\alpha\, {\bf w}_k -{\mu}\triangle {\bf w}_k +\gamma ({\bf w}_{k-1} \cdot \nabla) {\bf w}_k +
(1-\gamma)\mbox{ curl }{\bf w}_{k-1} \times {\bf w}_k+ \nabla q_k={\bf f}\qquad \mbox{ in } \,\,\qquad \Omega,
\label{eq:10}\\
\quad\qquad\qquad\qquad\qquad\qquad\qquad\qquad\qquad\qquad\qquad\qquad\qquad\,\mbox{ div } {\bf w}_k ={ 0} \qquad \mbox{ in } \qquad \,\,\Omega,
\label{eq:11}\\
\quad\qquad\qquad\qquad\qquad\qquad\qquad\qquad\qquad\qquad\qquad\qquad\qquad\quad\,\,\,\,\,\, {\bf w}_k =\bold{\sigma} \,\, \qquad \mbox{on} \qquad \partial \Omega
\label{eq:12}
\end{gather}
$ k = 1,2, \ldots, $ converges to the exact solution of the nonlinear problem
(7)-(9) for $ k \rightarrow \infty $ ($\gamma=1$ \cite{Elman} and $\gamma=0$ \cite{Benzi})
   for an arbitrary given initial approximation $ {\bf w}_0 $ satisfying the conditions
   $ \mbox {div} {\bf w}_0 = {0} \quad \mbox{in} \quad \Omega \, \, $ and
   $ \, \, {\bf w}_0 = \bold{\sigma} \quad \mbox {on} \quad \partial \Omega. $
Note that there is no need to set the initial approximation $ q_0 $ for pressure.

Thus, in order to find a solution to the original nonlinear problem (7)-(9), it is necessary for each
   $ k = 1,2, \ldots, $ in (10)-(12) find a solution to the linearized problem of the following form:
   find the velocity vector function
   $ {\bf u} = {\bf u} ({\bf x}) = (u_1 ({\bf x}), u_2 ({\bf x})) $ and
pressure $ p = p ({\bf x}) $, satisfying the system of equations and boundary conditions
\begin{gather}
\alpha\, {\bf u} -{\mu}\triangle {\bf u} +\gamma ({\bf g}_{k-1} \cdot \nabla) {\bf u}+
(1-\gamma)\mbox{ curl }{\bf z}_{k-1} \times {\bf u}+
\nabla p={\bf f}\qquad \mbox{ in } \,\,\qquad \Omega,
\label{eq:13}\\
\,\,\quad\qquad\qquad\qquad\qquad\qquad\qquad\qquad\qquad\qquad\qquad\,\,\,\,\mbox{ div } {\bf u} ={ 0} \qquad \mbox{ in } \qquad \,\,\Omega,
\label{eq:14}\\
\,\,\,\quad\qquad\qquad\qquad\qquad\qquad\qquad\qquad\qquad\qquad\qquad\qquad\,\,\, {\bf u} =\bold{\sigma} \qquad \,\, \mbox{on} \qquad \partial \Omega.
\label{eq:15}
\end{gather}
Here $ {\bf g}_{k-1} $ and ${\bf z}_{k-1 } $ are the approximation of the velocity vector function from the previous Picard iteration relevant case $\gamma=1$ and $\gamma=0$ respectively,
  $ {\bf g}_{k-1} \in L_{\infty} (\Omega) \times L_{\infty} (\Omega), $
$ \mbox{div } {\bf g}_{k-1} = 0 $ in $ \Omega $ and $ \mbox{ curl } {\bf z}_{k-1} \in L_{\infty} (\Omega), $
$ \mbox{div } {\bf z}_{k-1} = 0 $ in $ \Omega $.

The peculiarity of studying a sequence of systems of the form (\ref{eq:13})-(\ref{eq:15}) ($ k = 1,2, \ldots $),
and therefore nonlinear
problem (\ref{eq:7})-(\ref{eq:9}) will consist in the fact that the domain $ \Omega $ is
a non-convex polygon with a reentrant corner $ \omega, $ that is, a corner greater than 180 degrees,
 whose vertex coincides with the origin $ {\cal O} = (0, 0) $.

It is well known (see, for example, \cite{Girault}) that approximately solving
a similar system (\ref{eq:13})-(\ref{eq:15}) by classical finite-difference
or by finite element methods in the domain with the reentrant corner $ \omega, $
the error arising in the vicinity of its vertex propagates into the inner part of the computational domain
 $  \Omega, $ where the solution of the problem has sufficient smoothness.
 And as a consequence, the result (for the velocity vector function)  is poor from the point of view of convergence,
namely, if
$ {\bf u}_h $ is an approximate solution to the problem (\ref{eq:13})-(\ref{eq:15}), then the estimate
\begin{gather}
\|\nabla({\bf u}-{\bf u}_h)\|_{{\bf L_2}(\Omega_h)}+h^{-\lambda} \|{\bf u}-{\bf u}_h\|_{{\bf L_2}(\Omega_h)}
\leq C h^{\lambda},
\label{eq:16}
\end{gather}
where $ \lambda <1 $, $ C $ is a positive constant independent of $ {\bf u} ({\bf x}) $ and $ h. $
Recall that for convex polygonal domains the order $ \lambda $ in the estimate (\ref {eq:16})
is equal to one, and for nonconvex polygonal domains $ \Omega $, the order of convergence $ \lambda $
becomes much smaller as the reentrant  corner $ \omega $ increases. For example,
if  $ \omega $ is $202.5$ degrees, then the exponent $ \lambda $ is approximately $ 0.8 $, and if $ \omega = 270^{\circ}, $
then $ \lambda $ is already --- $ 0.54, $ that is with such a value of the corner,
the order of convergence is almost halved with respect to $ \, h $,
in comparison with the case of a convex domain.
Consequently, to achieve a given accuracy, it is necessary to use a grid with a much smaller value of its step $ h $ and,
as a consequence, very significant time and resource costs are required.

In the presented paper, we will define the notion of $ R_{\nu}$-generalized solution in weighted sets.
Let's build a weighting method
 finite elements and on its basis we obtain a sequence of approximate solutions converging to the solution of the
nonlinear problem (\ref{eq:7})-(\ref{eq:9}) with the first order with respect to $ \, h \, $
 for domains
 with different values
    of the reentrant corner $ \omega $
   greater than $ \pi $ for $\gamma=1$ and $\gamma=0$.
 Thus, we will experimentally establish the independence of the order of convergence of the approximate solution to the
exact one of the problem (\ref{eq:7})-(\ref{eq:9}) from the value of the reentrant corner.

 In order to define the $ R_{\nu}$-generalized solution of a sequence of problems of the form (\ref{eq:13})-(\ref{eq:15}),
 we introduce the necessary spaces and sets of functions. By $ \Omega^{\delta} $ we denote the intersection of the disk
of radius $ \delta, \delta> 0, $
 centered at the origin $ {\cal O} = (0, 0) $ with the closure of the domain $ \Omega. $
Define the function $ \rho ({\bf x}) $ in $ \bar {\Omega}, $
 which we will call weighted function, satisfying the following conditions: $ \rho ({\bf x}) = \sqrt {x_1 ^ 2 + x_2 ^ 2}, $
 if $ {\bf x} \in \Omega^{\delta} $ and $ \rho ({\bf x}) = \delta $ otherwise.
Let $ D^l v ({\bf x}) = \frac {\partial^{|l|} v ({\bf x})}
 {\partial x_1^{l_1} \partial x_2^{l_2}} $, $ l = (l_1, l_2), $ $ |l| = l_1 + l_2, \, l_i $ - non-negative integers
$ i \in \{1, 2 \} $,
 $ d {\bf x} = d x_1 d x_2 $.

We denote by $L_{2, \beta} (\Omega),
W^1_{2, \beta} (\Omega)
$
the space of functions $ v ({\bf x}) $ such that
\begin{gather}
\|v\|_{L_{2,\beta}(\Omega)}=\Bigr(\int\limits_{\Omega}{\rho^{2\beta}({\bf x})\,
v^2({\bf x})  d{\bf x}}\Bigl)^{1/2},
\label{eq:17}\\
\|v\|_{W^1_{2,\beta}(\Omega)}=
\Bigr(\|v\|^2_{L_{2,\beta}(\Omega)}+\sum\limits_{|l|=1}{\int\limits_{\Omega}{\rho^{2\beta}({\bf x})
|D^l v({\bf x})|^2  d{\bf x}}}\Bigl)^{1/2}
\label{eq:18}
\end{gather}
respectively.

Let us define the following conditions for the functions $ v({\bf x}) $
\begin{gather}
0 < C_1\leq \|v\|_{L_{2,\beta}(\Omega\setminus\Omega^{\delta})},
\label{eq:19}\\
|D^l v({\bf x})| \leq C_2 \Bigl(\delta \rho^{-1}({\bf x})\Bigr)^{\beta+l-\varepsilon} {\bf x} \in \Omega^{\delta},
\label{eq:20}
\end{gather}
where $ C_2 $ is a positive constant independent of $ v ({\bf x}) $ and $ l, l \ in \{0, 1 \}, $
$ \varepsilon $ is a small positive parameter that does not depends on $ \delta, \beta, l $ and $ v ({\bf x}). $

By $ L_{2,  \beta} (\Omega,  \delta) $
we denote the set of functions $ v ({\bf x}) $ from the space $ L_{2, \beta} (\Omega) $
satisfying the conditions (\ref{eq:19}) and (20) ($ l = 0 $) with bounded norm (17).
 Define the subset $ Y_{\beta} (\Omega, \delta) =
\{
v({\bf x}) \in L_{2, \beta} (\Omega, \delta):
\| \rho^{\beta} \, v \|_{L_1 (\Omega)} = 0 \} $ with limited norm (17).

Let $ X_{\beta} (\Omega, \delta) $ be the set of functions $ v ({\bf x}) $ from the space
$ W^1_{2, \beta} (\Omega) $ satisfying the conditions (19) and (20) with bounded norm
(18). Through $ X^0_{\beta} (\Omega, \delta) $
we denote the closure in the norm (18) of the set of infinitely differentiable functions
with compact support in $ \Omega $ satisfying the conditions (19), (20). Next,
define the set $ Z_{\beta} (\partial \Omega, \delta) = \{s ({\bf x}): \mbox{ such that there is } S({\bf x}) \in
X_{\beta} (\Omega, \delta) \mbox{ and } S({\bf x}) = s({\bf x}) \mbox{ on } \partial \Omega \} $
with bounded norm
\begin{gather}
\|s\|_{Z_{\beta}(\partial\Omega, \delta)}=
\inf\limits_{S({\bf x})=s({\bf x}) \mbox{ on } \partial \Omega}{\|S\|_{W^1_{2,\beta}(\Omega)}}.
\label{eq:21}
\end{gather}

We will use bold type to highlight the sets of vector functions
${\bf L}_{2,\beta}(\Omega,\delta)=\{ {\bf v}({\bf x})= (v_1({\bf x}), v_2({\bf x})): v_i({\bf x})\in
 {L}_{2,\beta}(\Omega,\delta)\}$
with bounded norm
$\|{\bf v}\|_{{\bf L}_{2,\beta}(\Omega)}=
 \Bigl(\|v_1\|^2_{L_{2,\beta}(\Omega)}+\|v_2\|^2_{L_{2,\beta}(\Omega)}\Bigr)^{1/2}.$
The sets of vector functions
${\bf X}_{\beta}(\Omega, \delta)$ (${\bf X}^0_{\beta}(\Omega, \delta)$)
 with vector norm (18) and ${\bf Z}_{\beta}(\partial\Omega, \delta)$  with vector norm
 (21)
are introduced similarly.

We assume that the vector functions of the right-hand sides (13) and (15) satisfy the conditions
\begin{gather}
  {\bf f}({\bf x})\in {\bf L}_{2,\beta}(\Omega,\delta),\qquad {\bf \sigma}({\bf x})\in
{\bf  Z}_{\beta}(\partial \Omega,\delta) \qquad \beta>0.
\label{eq:22}
\end{gather}

Define bilinear and linear forms
\begin{displaymath}
a^{\gamma}_{(k)}({\bf y}, {\bf v})=
\int\limits_{\Omega}{
\Bigl[
\alpha {\bf y}\cdot(\rho^{2\nu}{\bf v})+ \mu
\nabla {\bf y}: \nabla(\rho^{2\nu} {\bf v})+ \gamma({\bf g}_{k-1} \cdot \nabla)
{\bf y}\cdot (\rho^{2\nu}{\bf v})  +(1-\gamma) (\mbox{curl} {\bf z}_{k-1} \times {\bf y})\cdot (\rho^{2\nu}{\bf v})\Bigr] d {\bf x}},
\end{displaymath}
\begin{displaymath}
b_1({\bf v}, s)=-\int\limits_{\Omega}{s\, \mbox{ div }(\rho^{2\nu}{\bf v}) d {\bf  x}}=b({\bf v}, s)+c({\bf v}, s),
\end{displaymath}
\begin{displaymath}
b_2({\bf v}, s)=-\int\limits_{\Omega}{(\rho^{2\nu}\, s)\,  \mbox{ div }{\bf v}\,  d {\bf x}}=b({\bf v}, s)-c({\bf v}, s),
\end{displaymath}
where
\begin{gather}
 b({\bf v}, s)=-\int\limits_{\Omega}{(\rho^{\nu}\,s) \, \mbox{ div }(\rho^{\nu}{\bf v}) d {\bf  x}}, \qquad
 c({\bf v}, s)=-\int\limits_{\Omega}{(\rho^{\nu} s) ({\bf v}\cdot \nabla \rho^{\nu} )d {\bf x}},
\label{eq:23}
\end{gather}
\begin{displaymath}
l({ \bf v})=\int\limits_{\Omega}{ {\bf f}\cdot ( \rho^{2\nu} {\bf v}) d {\bf x}}
\end{displaymath}
and introduce the concept of an $ R_{\nu} $-generalized solution (13)-(15)
 for $ k = 1, 2, \ldots. $

{\bf Definition 1}. {\it An $ R_{\nu} $-generalized solution of the problem (13)-(15)
is a pair of functions
$ ({\bf u}_{\nu}, \, p_{\nu}) \in {\bf X}_{\nu} (\Omega, \delta) \times Y_{\nu} (\Omega, \delta), $
$ {\bf u}_{\nu} $ satisfies almost everywhere the condition (15) on $ \partial \Omega, $
such that for all pairs
  $ ({\bf v}, g) \in {\bf X}^0_{\nu} (\Omega, \delta) \times Y_{\nu} (\Omega, \delta) $
  the integral identities
\begin{gather}
a^{\gamma}_{(k)}({\bf u}_{\nu}, {\bf v})+b_1({\bf v}, p_{\nu})=l({\bf v}),
\label{eq:24}\\
 b_2({\bf u}_{\nu}, g)=0\qquad\qquad\qquad\,\,\,\,\,
\label{eq:25}
\end{gather}
 hold, where
$ {\bf f} $ and $ {\bf \sigma} $ are well-known vector functions in $ \Omega $ and on $ \partial \Omega $,
respectively, which satisfy the conditions (22), for $ \nu \geq \beta $}.

{\bf Remark 2}. {\it In \cite{Ruk1989} it is proved that the function $ \rho^{\nu} s \in L_{2,0} (\Omega, \delta) $
if and only if the function $ s \in L_{2, \nu} (\Omega, \delta) $ similar to $ \rho^{\nu} s \in X_{0}(\Omega, \delta)
 \Leftrightarrow s \in X_{\nu} (\Omega, \delta) $ and $ \rho^{\nu} s \in X^0_ {0} (\Omega, \delta) \Leftrightarrow s
 \in X^0_{\nu} (\Omega, \delta)$. Therefore, the bilinear form $ b(\cdot, \cdot) $ defined in (23)
 is the classical bilinear form
for saddle problems (see, for example, \cite{Brezzi})}.

{\bf Remark 3}. {\it Besides the fact that the $ R_{\nu} $-generalized solution is defined in weighted sets.
 It is easy to see that, due to the definition of the weight function $ \rho({\bf x}) $, the bilinear form
$ c(\cdot, \cdot) $ (see (23)):
  $ c({\bf v}, s) = - \int \limits_{\Omega}{(\rho^{\nu} s) ({\bf v} \cdot \nabla \rho^{\nu}) d {\bf x}} = -
\int \limits_{\Omega^{\delta}} {(\rho^{\nu} s) ({\bf v} \cdot \nabla \rho^{\nu} ) d {\bf x}} $
and in the general case is not equal to zero. Thus, $ b_1({\bf v}, s) \neq b_2 ({\bf v}, s) $
and individually they are not equal to $ b({\bf v}, s). $}

\paragraph{2. Construction of an approximate method for the problem.}

Let's perform a quasi-uniform triangulation $ \Lambda_h $ of the domain $ \bar {\Omega} $ \cite{Ciarlet}.
We split
domain into triangles $ \, K_j \, $ with side of order $ h, $ which we call basic elements.
Next, connect the vertices of the basic element with its center of gravity.
 Thus, we divide the element $ \, K_j \, $ into three triangles $ L_{j_s} $, which we
call finite elements, their set over all $ j $ forms a triangulation of $ T_h $ and
$ \Omega_h = \bigcup \limits_{L_{j_s} \in T_h} {L_{j_s}} $.

Let us define the sets of nodes for the components of the velocity and  pressure,  and introduce for each of them
finite-dimensional spaces: \\
{\it 1. For the components of the velocity.} As the approximating nodes we take the vertices of the finite elements
and the midpoints of the sides $ L_{j_s} $. Using the selected nodes, we define the basis functions
$ \theta_{(n)} ({\bf x}), {\bf x} \in L_{j_s} $
which are polynomials of degree two. Moreover, we assume that the coinciding nodes of neighboring elements are common.
Therefore, the functions in the presented basis will be continuous in $ \bar {\Omega} $.
 Let us denote the space of them by $ V_h, $ and the space of vector functions by $ {\bf V}_h = V_h \times V_h $.
  $ {\bf V}_h $ coincides with the space $ {\bf P}_2 $ of
the well-known Taylor-Hood pair $ {\bf P}_2 $ - $ P_1 $ (see \cite{Brezzi}). \\

{\it 2. For the pressure.} We take the vertices of the finite elements $ L_{j_s} $ as approximating nodes.
Using these nodes, we define the basis functions $ \chi_{(l)} ({\bf x}), {\bf x} \in L_{j_s}, $
which are polynomials of the first degree. We will assume that the coinciding nodes of two (or more) neighboring
 finite elements are different nodes,
functions represented in such a basis will be discontinuous and only belong to the space $ L_2 (\Omega) $.
We denote their totality by $ W_h. $ Note that the space $ W_h $ differs from the continuous space $ P_1, $ of the Taylor-Hood pais $ {\bf P}_2 $ - $ P_1, $
and its dimension is substantially greater than the dimension of the $ P_1 $ one.
The payment for this inconvenience is the fact that for any vector function
$ {\bf v}_h \in {\bf V}_h $ there is such a function $ s_h \in W_h, $ that $ s_h = \mbox{div} { \bf v}_h $
and from the condition $ \int \limits_{\Omega} {(\rho^{2 \nu} \, s_h) \, \mbox{div} {\bf v}_h \, d {\bf x}} = 0 $
follows
$ \| \mbox{div} {\bf v}_h \|_{L_{2, \nu} (\Omega)} = 0. $ In addition,  functions of the pair of spaces
$ {\bf V}_h $ -- $ W_h $ satisfy $ inf-sup $ condition \cite{Qin}.

Further, we improve the basis functions $ \theta_{(n)} ({\bf x}) $ and $ \chi_{(l)} ({\bf x}) $
of the introduced spaces. For this purpose, we multiply them by the weight function $ \rho ({\bf x}) $ in some degree
\begin{gather}
\varphi_{(n)} ({\bf x}) = \rho^{\nu^{\ast}} ({\bf x}) \, \theta_{(n)}({\bf x}) \qquad \qquad \mbox {and} \qquad \qquad
\psi_{(l)} ({\bf x}) = \rho^{\mu^{\ast}} ({\bf x}) \, \chi_{(l)} ({\bf x}) ,
\label {eq:26}
\end{gather}
where $ \nu^{\ast} $ and $ \mu^{\ast} $ are real numbers (free parameters of the method).
Looking ahead, note that the optimal error of the method is achieved when we take them equal each other.

Next, we define the finite-dimensional space of functions $ X_h $ for the components of the velocity
$ {\bf u}_{h, \nu} ({\bf x}) = (u_{h, \nu, 1} ({\bf x}), u_{h, \nu, 2} ({ \bf x})) $ such that
\begin{gather}
u_{h,\nu,1}({\bf x})=\sum\limits_{i=1}^{dim\, V_h}{\hat{c}_{i}\, \varphi_{(i)}({\bf x})},\qquad \qquad \qquad
u_{h,\nu,2}({\bf x})=\sum\limits_{i=1}^{dim\, V_h}{\hat{d}_{i}\, \varphi_{(i)}({\bf x})},
\label{eq:27}
\end{gather}
where $ \varphi_{(i)}({\bf x}) $ are basis functions  (see (26)),
  $ \hat{c}_i = \rho^{- \nu^{\ast}}(M_i) \, c_{i} $ and $ \hat{d}_i = \rho^{- \nu^{\ast}} (M_i) \,
 d_{i} $ are coefficients, $ M_i $ are nodes of approximation of the space $ V_h $ in $ \Omega $.
 We denote by $ \, X^0_h \, $ the subspace
$ X_h: $ $ X^0_h = \{v_{h} \in X_h: v_{h} (M_i) = 0, M_i \in \partial \Omega \}. $

Then, we define the finite-dimensional space $ Y_h $ for the pressure
$ p_{h, \nu} ({\bf x}) $ such that
\begin{gather}
p_{h,\nu}({\bf x})=\sum\limits_{j=1}^{dim\, W_h}{\hat{e}_{j}\, \psi_{(j)}({\bf x})},
\label{eq:28}
\end{gather}
where $ \psi_{(j)} ({\bf x}) $ are basis functions  (see (26)),
 $ \hat{e}_{j} = \rho^{- \mu^{\ast}}(N_j) e_{j} $ are coefficients, $ N_j $ are approximation nodes of the space $ W_h $.
  Moreover, $ {\bf X}_h =X_h \times X_h \subset{\bf X}_{\nu} (\Omega, \delta), {\bf X}^0_h = X^0_h \times X^0_h \subset
 {\bf X}^0_{\nu} (\Omega, \delta) $ and $ Y_h \subset Y_{\nu} (\Omega, \delta) $.

Now we define an approximate $ R_{\nu} $ - generalized solution of the problem (13)-(15),
$ k =1,2, \ldots.$

{\bf Definition 2}. {\it An approximate $ R_{\nu} $-generalized solution of the problem (13)-(15)
is a pair of functions
$ ({\bf u}_{h, \nu}, \, p_{h, \nu}) \in {\bf X}_{h} (\Omega, \delta) \times Y_{h} ( \Omega, \delta), $
$ {\bf u}_{h, \nu} $ satisfies the condition (15) at nodes on $ \partial \Omega, $
  such that for all pairs
  $ ({\bf v}_h, g_h) \in {\bf X}^0_{h} \times Y_{h} $
  the integral identities
\begin{gather}
a^{\gamma}_{(k)}({\bf u}_{h,\nu}, {\bf v}_h)+b_1({\bf v}_h, p_{h,\nu})=l({\bf v}_h),
\label{eq:29}\\
 b_2({\bf u}_{h,\nu}, g_h)=0\qquad\qquad\quad\,\,\,\qquad\quad
\label{eq:30}
\end{gather}
 hold, where
$ {\bf f} $ and $ {\bf \sigma} $ obey the conditions (22), $ \nu \geq \gamma $}.

In order to find an approximate solution of the problem (13)-(15)
 in the statement (29), (30), it is necessary solve a system
of linear algebraic equations of the following form:
\begin{gather}
 {\bf A}^{\gamma}_{(k)} {\bf y}+{\bf B}_1 {\bf z}={\bf F},
\label{eq:31}\\
{\bf B}_2^T {\bf y}={\bf 0},\qquad\qquad
\label{eq:32}
\end{gather}
here $ {\bf y}=(\hat{c}_1, \hat{c}_2, \ldots, \hat{c}_{dim \, V_h}, \hat{d}_1, \hat{d}_2, \ldots , \hat{d}_{dim \, V_h})^T,
 \, {\bf z} = (\hat{e}_1, \hat{e}_2, \ldots, \hat{e}_{ dim \, W_h})^T $ and
$ {\bf F} $ is a column vector of dimension $ 2 \times dim V_h $ of values of the linear form $ l (\varphi_{(i)}) $.

{\bf Remark 4}. {\it Solving the system of equations (31), (32),
we find the unknowns $ \hat{c}_i, \hat {d}_i $ and $ \hat{e}_j . $ Then, using the reverse formulas (see (27),
 (28)
 $ c_i = \rho^{\nu^{\ast}} (M_i) \, \hat{c}_i, d_i = \rho^{\nu^{\ast}} (M_i) \, \hat{d}_i $ and
$ e_{j} = \rho^{\mu^{\ast}} (N_j) \hat{e}_j $ we calculate the values of the velocity and the pressure at the approximation nodes $ M_i $ and $ N_j $. }

To find the solution $ ({\bf y}, {\bf z}) $ of the system of equations (31), (32), for a fixed $ k $, we
apply the converging Uzawa \cite{Bramble} incomplete procedure with block preconditioning of its matrix.
We represent the procedure in the form of Algorithm 1, the stages of it are implemented using auxiliary Algorithms 2 and 3.

{\bf Algorithm 1}. \\
{\it Stage 1}. Let an initial approximation $ ({\bf \chi}^0, \, {\bf \xi}^0) $
 satisfying the equations of the system (31), (32) be given. \\
{\it Stage 2}.
a) Calculate the vector
 \begin {gather}
\tilde {\bf r}^l = {\bf F} - {\bf A}^{\gamma}_{(k)}{\bf \chi}^l - {\bf B}_1{\bf \xi}^l .
\label{eq:33}
 \end{gather}
b) Find the following approximation
 \begin {gather}
{\bf \chi}^{l+1}={\bf\chi}^l+{\bf r}^l,
\label{eq:34}
 \end{gather}
 applying the generalized minimum residual method (GMRES-method) with the preconditioner
$ \hat{\bf A}^{\gamma}_{(k)} $ to the matrix
 $ {\bf A}^{\gamma}_{(k)} $. The calculation of the vector $ {\bf r}^l $ in (34) is described in Algorithm 2. \\
{\it Stage 3}.
a) Calculate the vector
 \begin{gather}
\tilde{\bf d}^l = {\bf B}_2^T {\bf \chi}^l.
\label{eq:35}
 \end{gather}
b) Find the following approximation
 \begin{gather}
{\bf\xi}^{l+1} = {\bf\xi}^l+{\bf d}^l,
\label{eq:36}
 \end{gather}
 applying an internal iterative procedure, at each iteration of which the GMRES method was used. A
  preconditioner $ \hat {\bf S}^{\gamma}_{(k)} $ to $ {\bf S}^{\gamma}_{(k)}
= {\bf B}_2^T \, ({\bf A}^{\gamma}_{(k)})^{-1} \, {\bf B}_1 $ using an auxiliary matrix
$ \tilde {\bf S}^{\gamma,0}_{(k)} $ is constructed. The calculation of the vector $ {\bf d}^l $ in (36)
 is described in Algorithm 3. \\
 {\it Step 4}.
 Check if the stop condition is met. If satisfied, then we have found a solution $ ({\bf \chi}^{L_k}, \,
 {\bf \xi}^{L_k}). $ Otherwise, go to Stage 2.

{\bf Algorithm 2}. \\
 The preconditioning matrix $ \hat{\bf A}^{\gamma}_{(k)} $ to $ {\bf A}^{\gamma}_{(k)} $ is obtained
by incomplete decomposition into the lower triangular $ \hat{\bf L}^{\gamma}_{(k)} $ and upper triangular
$ \hat {\bf U}^{\gamma}_{(k)} $ matrices (see, for example, \cite{Saad}): $ \hat {\bf A}^{\gamma}_{(k)} = \hat {\bf L}^{\gamma}_{(k)}
\cdot \hat {\bf U}^{\gamma}_{(k)}, \, \, (\hat{\bf A}^{\gamma}_{(k)})^{-1} =
  (\hat{\bf U}^{\gamma}_{(k)})^{-1} \cdot (\hat {\bf L}^{\gamma}_{(k)})^{-1} $.

 Let $ {\bf q} = (\hat{\bf A}^{\gamma}_{(k)})^{-1} \tilde{\bf r}^l, $ where the vector $ \tilde {\bf r}^l $
computed in (33),
 then the Arnoldi procedure (\cite{Saad}) generates an orthogonal basis with a left preconditioned Krylov subspace
of dimension $ s $ with minimal residual
 $ {\bf r}^l = \sum\limits_{i = 0}^{s-1} {\beta_i \Bigl ((\hat{A}^{\gamma}_{(k)})^{-1} \, A^{\gamma}_{(k)} \Bigr)^i} {\bf q},\,\,
 \beta_i \in {\bf R}, s = 5. $

{\bf Algorithm 3}. \\
  Suppose we need to calculate the vector $ (\hat {\bf S}^{\gamma}_{(k)})^{-1} \tilde {\bf d}^l, $
  but finding the inverse of $ \hat {\bf S}^{\gamma}_{(k)} $ is resource-intensive. It is much easier to multiply the auxiliary
  the matrix $ {\bf S}^{\gamma,0}_{(k)} $ to it by a vector. The components of which are determined using the integrals
\begin{displaymath}
\Bigl({\bf S}^{\gamma,0}_{(k)} \Bigr)_{ij} = \frac{1}{\mu} \int \limits_{\Delta_{ij}} {\rho^{2 \, \nu + 2\, \mu^{\ast}}
 \chi_{(i)}({\bf x}) \, \chi_{(j)} ({\bf x}) d {\bf x}},
\end{displaymath}
where $ \Delta_{ij} = \mbox{supp} \{\chi_{(i)} ({\bf x}) \} \cap \mbox{supp} \{\chi_{(j)} ({\bf x}) \} \cap \Omega $
and
$ \chi_{(i)} ({\bf x}), \chi_{(j)} ({\bf x}) $ are basis functions of the space $ W_h$. By $ \tilde {\bf S}^{\gamma}_{(k)} $
we denote the diagonal matrix, elements of which are of the form
$ \Bigl(\tilde{\bf S}^{\gamma}_{(k)} \Bigr)_{ii} = \sum\limits_{j = 1}^{dim W_h} {\Bigl({\bf S}^{\gamma,0}_{(k)} \Bigr)_{ij}} $.
We organize an iterative procedure, at each iteration of which, using the GMRES method, we find the vector
by analogy with $ {\bf r}^l $ of Algorithm 2: \\
1. Let $ {\bf e}^0 $ be the initial approximation, $ {\bf e}^0 = (0,0, \ldots, 0)^T. $ \\
2. For $ n = 0,1, \ldots, N-1, $ calculate
\begin{displaymath}
{\bf e}^{n+1} = {\bf e}^{n} + (\tilde {\bf S}^{\gamma}_{(k)})^{-1} \,(\tilde{\bf d}^l - {\bf S}_{(k)}^{\gamma,0} {\bf e}^{n}).
\end{displaymath}
3. The vector $ {\bf e}^{N} $ is the required vector $ {\bf d}^l $ in (36).

Using Algorithm 4 below, we
find the values of the approximate solution of the nonlinear problem (7)-(9)
at the grid nodes and, therefore,
the velocity and the pressure functions, which we denote by $ {\bf w}_{h, \nu} $
and $ q_{h, \nu} $ respectively.

{\bf Algorithm 4}.\\
{\it Stage 1}.
Let the vector $ {\bf y}_0 $ be an initial approximation satisfying the conditions (7) and (8)
at the nodes $ M_i $ in $ \Omega $ and at $ \partial \Omega $, respectively. Without loss of generality, we assume
that $ {\bf y}_0 = {\bf 0} $ and $ {\bf z}_0 = {\bf 0} $ in the nodes $ M_i $ and $ N_j $ in $ \Omega $, respectively. \\
{\it Stage 2}.
We implement Picard's iterative procedure, $ k = 1,2, \ldots: $ \\
1. Denote by $ (\chi^0, \, \xi^0) $ a pair of vectors $ ({\bf y}_{k-1}, \, {\bf z}_{k-1})$. \\
2. Let's execute Algorithm 1. \\
3. Denote by $ ({\bf y}_{k}, \, {\bf z}_{k}), $ the pair of vectors obtained as a result of applying Algorithm 1
$ (\chi^{L_k}, \, \xi^{L_k})$. \\
{\it Stage 3}.
If the stopping condition is met, then the desired solution has been found. Otherwise, go to stage 2.

{\bf Remark 5}.

{\it If $ \nu = \nu^{\ast} = \mu^{\ast} = 0 $, then we have a sequence of classical
approximate solutions (29), (30). Using Algorithm 4 we get a resulting pair $ ({\bf w}_h, \, q_h) $.}

\paragraph{3. The results of numerical simulation.}

In this section, we present the results of numerical experiments for an approximate solution for the weighted
FEM of the nonlinear problem
(7)-(9) in convective form for various non-convex polygonal domains $ \Omega_m $ of the form
$ (-1,\, 1) \times (-1,\, 1) \setminus \bar{J}_m $ with one reentrant corner $ \omega $ taking the values $ \omega_m $
greater than $ \pi $ on the boundary with the vertex at the origin.
 Let us compare the order of convergence of the method with the corresponding classical FEM
$ (\nu = \nu^{\ast} = \mu^{\ast} = 0) $. The corner $ \omega $ in the tests will take the values      $\omega_m=\Bigl(1+\frac{1}{2^m}\Bigr)\,\pi$ and the domains $ \bar{J}_m, m = 1,2,3 $
for each of them have the form $ \bar{J}_1 = \{(x_1, \, x_2): -1<x_1<0, -1<x_2<0 \},
\bar{J}_2 = \{(x_1, \, x_2): -1<x_1<0, -1<x_2<x_1 \}, \bar{J}_3 = \{(x_1, \, x_2): -1<x_1<0, -1<x_2<\frac{x_1}{2} \}. $

In tests, the exact solution of the nonlinear problem (7)-(9) in polar coordinates
$ (r, \varphi) $ depends on the reentrant corner  $ \omega_m $ has the form:
\begin{gather}
w_1(r, \varphi)=r^{\lambda_m}\, \theta_1(\varphi),\quad w_2(r, \varphi)=r^{\lambda_m}\, \theta_2(\varphi),\quad
q(r, \varphi)=r^{\lambda_m-1}\, \theta_3(\varphi),\quad
\label{eq:37}
 \end{gather}
 where
\begin{gather}
\left (\begin{array}{c}\displaystyle
{\theta_1(\varphi)}
\\\displaystyle {\theta_2(\varphi)}
\end{array}\right )
=
\left (\begin{array}{cc}\displaystyle (\lambda_m+1) \Lambda_m(\varphi) & \Lambda_m'(\varphi)
\\\displaystyle \Lambda_m'(\varphi) & - (\lambda_m+1) \Lambda_m(\varphi)
\end{array}\right)
\left (\begin{array}{c}\displaystyle {\sin \varphi}
\\\displaystyle {\cos \varphi}
\end{array}\right )\label{eq:38}\\
\theta_3(\varphi)=\frac{1}{\lambda_m-1}\Bigr( (\lambda_m+1)^2 \, \Lambda_m'(\varphi)+ \Lambda_m'''(\varphi)\Bigr), m=1,2,3.
\label{eq:39}
\end{gather}

The function $ \Lambda_m (\varphi) $ in the representations (38), (39) has the form
\begin{displaymath}
\Lambda_m(\varphi)=\cos ((\lambda_m-1)\varphi)-\cos ((\lambda_m+1)\varphi) + \cos (\lambda_m\cdot \omega_m)
\Bigl(\frac{\sin((\lambda_m+1)\varphi)}{\lambda_m+1}-  \frac{\sin((\lambda_m-1)\varphi)}{\lambda_m-1}\Bigr).
\end{displaymath}
$ \Lambda_m'(\varphi) $ and $ \Lambda_m''' (\varphi) $ ---
the first and third derivatives of the function
$ \Lambda_m(\varphi) $ by variable $ \varphi $, respectively; $ \lambda_m $ --- the smallest
positive number that is a solution to the following equation $ \lambda \sin \omega_m =
-\sin (\lambda \, \omega_m). $ Let us write out approximate values  for the corresponding $ \omega_m, m = 1,2,3 $
of the reentrant corner $ \omega: \omega_1 = 3 \pi/2 - \lambda_1 \approx 0.54448, \omega_2 = 5\pi/4 - \lambda_2
\approx 0.67358 $ and $ \omega_3 = 9 \pi/8 - \lambda_3 \approx 0.80077. $

It is well known that the solution $ ({\bf w}, q) $ is analytic in $ \bar {\Omega}_m \setminus {(0,0)}, $ but
$ {\bf w} \not \in {\bf W}^2_2 (\Omega_m) $ and $ q \not \in W_2^1 (\Omega_m). $ As a result, the error of the velocity vector function arising in the vicinity of the top of the reentrant corner, when using the classical FEM, propagates into the inner part of the calculated area. In this case, the order of convergence of the approximate solution to the exact solution of the problem significantly decreases with an increase in the value of the input corner $ \omega. $ Let $ h $ be the step of the quasi-uniform grid if $ \omega = 9 \pi/8, $ then
the order of convergence is $ {\cal O} (h^{0.8}), $ if $ \omega = 5 \pi/4, $ then
--- $ {\cal O} (h^{0.67}) $ and if $ \omega = 3 \pi/2, $ then
--- $ {\cal O} (h^{0.54}) $ in the norm of the space $ {\bf W}_2^1 (\Omega_m) $ Table 1, $\gamma=1$ ($\gamma=0$ the same).

The numerical method proposed by us contains the following free parameters: $ \delta $ ---
the value of the neighborhood of the corner $ \omega, $ the exponents $ \nu $ and $ \nu^{\ast} $
 in the definition of the $ R_{\nu} $ -generalized solution and the weighted functions in the basis, respectively.
The optimal choice of these parameters, determined as a result of a number of numerical experiments,
gives a significantly better result in relation to classical approaches, namely, the rate of convergence of the approximate
 $ R_{\nu}$-generalized solution to the exact one of the problem does not depend on the value of the reentrant corner
$ \omega $ and is equal to $ {\cal O} (h) $ in the norm of the space $ {\bf W}^1_{2, \nu} (\Omega) $.
In this case, the choice of the value $ \delta $ does not depend on the value of the $ \omega $ and is directly
proportional to the step of the grid $ h $ (see Tables 2-7).  Further, we denote by $ d^G_i(M_j) =|w_{h,i}(M_j)-w_i(M_j)| $ and
$ d^{R_{\nu}}_i(M_j) = |w_{h, \nu, i}(M_j)-w_i(M_j)| $
quantities the absolute error of the $i$-th component of the velocity of the generalized
$ (\nu = \nu^{\ast} = \mu^{\ast}=0) $ and $ R_{\nu}$-generalized solutions at the node $ M_j $ lying in
 $ \Omega $ respectively. Tables 8-9 show the results of the approaches
in terms of the number of nodes (fraction of their total number)
in which the error does not exceed the specified values
$\Xi_n, d^G_i(M_j)\leq \Xi_n$ and $d^{R_{\nu}}_i(M_j)\leq \Xi_n.$

\begin{table}[!t]
	\centering\small
	\caption{The errors of the standard approximate FEM  in the norm of the space $ {\bf W}^1_2(\Omega) $ (convective form).  }
	\vspace*{2mm}
	\begin{tabular}{|c|c|c|c|}
		\hline
		$h, \omega_m$ & $3\pi/2$ &  $5\pi/4$  & $9\pi/8$  \\\hline
		$10^{-2}$ & $2.063\cdot10^{-1}$ &   $1.004\cdot10^{-1}$ & $3.615\cdot10^{-2}$ \\\hline
		$5\cdot 10^{-3}$ & $1.429\cdot10^{-1}$  &  $6.309\cdot10^{-2}$  & $2.074\cdot10^{-2}$ \\\hline
		$2.5\cdot 10^{-3}$ & $9.806\cdot10^{-2}$ &  $3.957\cdot10^{-2}$  & $1.19\cdot10^{-2}$\\\hline
			\end{tabular}
\label{tabl_1}
\end{table}

\begin{table}[!ht]
\centering
\caption{The errors of  a weighted FEM   in the norm of the space ${\bf W}^1_{2, \nu} (\Omega) $
(convective form), $\omega=3 \pi/2.$
 }
\medskip
  \begin{tabular}{|c|c|c|c|c|c|c|}
    \hline
                & \multicolumn{3}{c|}{ $\delta=0.0124, \nu=1.8$ } & \multicolumn{3}{c|}{ $\delta=0.0127, \nu=2.0$}  \\
       \hline
$\nu^{\ast}$        & $h=10^{-2}$  & $h=5\cdot 10^{-3}$ & $h=2.5\cdot 10^{-3}$ & $h=10^{-2}$ & $h=5\cdot 10^{-3}$ & $h=2.5\cdot 10^{-3}$\\
                  \hline
$-0.5$ &$5.335e-5$&$2.656e-5$&$1.322e-5$&$2.269e-5$&$1.131e-5$&$5.627e-6$ \\
    \hline
    $-0.475$ &$4.929e-5$&$2.449e-5$&$1.218e-5$&$2.088e-5$&$1.042e-5$&$5.205e-6$ \\
    \hline
    $-0.45$ &$4.555e-5$&$2.274e-5$&$1.143e-5$&$1.924e-5$&$9.596e-6$&$4.768e-6$ \\
    \hline
$-0.425$ &$4.211e-5$&$2.107e-5$&$1.046e-5$&$1.798e-5$&$8.964e-6$&$4.494e-6$ \\
    \hline
$-0.4$ &$3.895e-5$&$1.95e-5$&$9.691e-6$&$1.714e-5$&$8.517e-6$&$4.233e-6$ \\
    \hline
$-0.375$ &$3.569e-5$&$1.794e-5$&$8.914e-6$&$1.616e-5$&$8.042e-6$&$4.005e-6$ \\
    \hline
$-0.35$ &$3.223e-5$&$1.618e-5$&$8.049e-6$&$1.516e-5$&$7.565e-6$&$3.756e-6$ \\
    \hline
$-0.325$ &$2.919e-5$&$1.449e-5$&$7.249e-6$&$1.31e-5$&$6.572e-6$&$3.269e-6$ \\
    \hline
$-0.3$ &$2.902e-5$&$4.44e-5$&$7.211e-6$&$1.244e-5$&$6.192e-6$&$3.077e-6$ \\
    \hline
$-0.275$ &$2.738e-5$&$1.371e-5$&$6.867e-6$&$1.203e-5$&$6.021e-6$&$2.986e-6$ \\
    \hline
$-0.25$ &$2.764e-5$&$1.379e-5$&$6.913e-6$&$1.169e-5$&$5.859e-6$&$2.912e-6$\\
    \hline
      \end{tabular}
\end{table}

\begin{table}[!ht]
\centering
\caption{The errors of  a weighted FEM   in the norm of the space ${\bf W}^1_{2, \nu} (\Omega) $
(convective form), $\omega=5 \pi/4.$
 }
\medskip
  \begin{tabular}{|c|c|c|c|c|c|c|}
    \hline
                & \multicolumn{3}{c|}{ $\delta=0.0124, \nu=1.8$ } & \multicolumn{3}{c|}{ $\delta=0.0127, \nu=2.0$}  \\
       \hline
$\nu^{\ast}$        & $h=10^{-2}$  & $h=5\cdot 10^{-3}$ & $h=2.5\cdot 10^{-3}$ & $h=10^{-2}$ & $h=5\cdot 10^{-3}$ & $h=2.5\cdot 10^{-3}$\\
                  \hline
$-0.4$ &$2.925e-5$&$1.451e-5$&$7.299e-6$&$1.362e-5$&$6.684e-6$&$3.317e-6$ \\
    \hline
    $-0.375$ &$2.696e-5$&$1.347e-5$&$6.893e-6$&$1.251e-5$&$6.307e-6$&$3.096e-6$ \\
    \hline
    $-0.35$ &$2.478e-5$&$1.233e-5$&$6.122e-6$&$1.145e-5$&$5.727e-6$&$2.871e-6$ \\
    \hline
$-0.325$ &$2.273e-5$&$1.139e-5$&$5.653e-6$&$1.044e-5$&$5.169e-6$&$2.553e-6$ \\
    \hline
$-0.3$ &$2.089e-5$&$1.054e-5$&$5.202e-6$&$9.631e-6$&$4.8e-6$&$2.418e-6$ \\
    \hline
$-0.275$ &$1.917e-5$&$9.576e-6$&$4.744e-6$&$8.654e-6$&$4.329e-6$&$2.148e-6$ \\
    \hline
$-0.25$ &$1.748e-5$&$8.807e-6$&$4.326e-6$&$7.851e-6$&$3.975e-6$&$1.957e-6$ \\
    \hline
$-0.225$ &$1.635e-5$&$8.149e-6$&$4.105e-6$&$7.121e-6$&$3.55e-6$&$1.791e-6$ \\
    \hline
$-0.2$ &$1.606e-5$&$7.972e-6$&$3.972e-6$&$6.865e-6$&$3.399e-6$&$1.676e-6$ \\
    \hline
$-0.175$ &$1.546e-5$&$7.783e-6$&$3.904e-6$&$6.554e-6$&$3.29e-6$&$1.662e-6$ \\
    \hline
      \end{tabular}
\end{table}

\begin{table}[!ht]
\centering
\caption{The errors of  a weighted FEM   in the norm of the space ${\bf W}^1_{2, \nu} (\Omega) $
(convective form), $\omega=9 \pi/8.$
 }
\medskip
  \begin{tabular}{|c|c|c|c|c|c|c|}
    \hline
                & \multicolumn{3}{c|}{ $\delta=0.0124, \nu=1.8$ } & \multicolumn{3}{c|}{ $\delta=0.0127, \nu=2.0$}  \\
       \hline
$\nu^{\ast}$        & $h=10^{-2}$  & $h=5\cdot 10^{-3}$ & $h=2.5\cdot 10^{-3}$ & $h=10^{-2}$ & $h=5\cdot 10^{-3}$ & $h=2.5\cdot 10^{-3}$\\
                  \hline
$-0.3$ &$1.659e-5$&$8.225e-6$&$4.124e-6$&$7.124e-6$&$3.542e-6$&$1.771e-6$ \\
    \hline
    $-0.275$ &$1.503e-5$&$7.476e-6$&$3.757e-6$&$6.474e-6$&$3.255e-6$&$1.619e-6$ \\
    \hline
    $-0.25$ &$1.364e-5$&$6.841e-6$&$3.399e-6$&$5.872e-6$&$2.918e-6$&$1.458e-6$ \\
    \hline
$-0.225$ &$1.222e-5$&$6.11e-6$&$3.046e-6$&$5.264e-6$&$2.623e-6$&$1.321e-6$ \\
    \hline
$-0.2$ &$1.094e-5$&$5.486e-6$&$2.725e-6$&$4.682e-6$&$2.325e-6$&$1.157e-6$ \\
    \hline
$-0.175$ &$9.898e-6$&$5.003e-6$&$2.482e-6$&$4.218e-6$&$2.108e-6$&$1.049e-6$ \\
    \hline
$-0.15$ &$8.722e-6$&$4.362e-6$&$2.17e-6$&$3.786e-6$&$1.884e-6$&$9.472e-7$ \\
    \hline
$-0.125$ &$8.092e-6$&$4.034e-6$&$2.024e-6$&$3.482e-6$&$1.727e-6$&$8.595e-7$ \\
    \hline
      \end{tabular}
\end{table}

\begin{table}[!ht]
\centering
\caption{The errors of  a weighted FEM   in the norm of the space ${\bf W}^1_{2, \nu} (\Omega) $
(rotation form), $\omega=3 \pi/2.$
 }
\medskip
  \begin{tabular}{|c|c|c|c|c|c|c|}
    \hline
                & \multicolumn{3}{c|}{ $\delta=0.0123, \nu=2.0$ } & \multicolumn{3}{c|}{ $\delta=0.0127, \nu=1.9$}  \\
       \hline
$\nu^{\ast}$        & $h=10^{-2}$  & $h=5\cdot 10^{-3}$ & $h=2.5\cdot 10^{-3}$ & $h=10^{-2}$ & $h=5\cdot 10^{-3}$ & $h=2.5\cdot 10^{-3}$\\
                  \hline
$-0.5$ &$2.155e-5$&$1.076e-5$&$5.34e-6$&$3.275e-5$&$1.631e-5$&$8.124e-6$ \\
    \hline
    $-0.475$ &$2.021e-5$&$1.012e-5$&$4.986e-6$&$3.005e-5$&$1.504e-5$&$7.547e-6$ \\
    \hline
    $-0.45$ &$1.857e-5$&$9.36e-6$&$4.614e-6$&$2.769e-5$&$1.399e-5$&$6.96e-6$ \\
    \hline
$-0.425$ &$1.712e-5$&$8.459e-6$&$4.248e-6$&$2.567e-5$&$1.285e-5$&$6.41e-6$ \\
    \hline
$-0.4$ &$1.598e-5$&$8.016e-6$&$3.983e-6$&$2.396e-5$&$1.2e-5$&$6.031e-6$ \\
    \hline
$-0.375$ &$1.44e-5$&$7.221e-6$&$3.636e-6$&$2.182e-5$&$1.099e-5$&$5.483e-6$ \\
    \hline
$-0.35$ &$1.296e-5$&$6.414e-6$&$3.219e-6$&$2.077e-5$&$1.041e-5$&$5.238e-6$ \\
    \hline
$-0.325$ &$1.249e-5$&$6.28e-6$&$3.124e-6$&$1.96e-5$&$9.807e-6$&$4.891e-6$ \\
    \hline
$-0.3$ &$1.212e-5$&$6.025e-6$&$2.982e-6$&$1.907e-5$&$9.483e-6$&$4.768e-6$ \\
    \hline
$-0.275$ &$1.171e-5$&$5.882e-6$&$2.895e-6$&$1.868e-5$&$9.326e-6$&$4.599e-6$ \\
    \hline
$-0.25$ &$1.185e-5$&$5.968e-6$&$3.005e-6$&$1.91e-5$&$9.529e-6$&$4.77e-6$\\
    \hline
      \end{tabular}
\end{table}

\begin{table}[!ht]
\centering
\caption{The errors of  a weighted FEM   in the norm of the space ${\bf W}^1_{2, \nu} (\Omega) $
(rotation form), $\omega=5 \pi/4.$
 }
\medskip
  \begin{tabular}{|c|c|c|c|c|c|c|}
    \hline
                & \multicolumn{3}{c|}{ $\delta=0.0123, \nu=1.9$ } & \multicolumn{3}{c|}{ $\delta=0.0127, \nu=1.8$}  \\
       \hline
$\nu^{\ast}$        & $h=10^{-2}$  & $h=5\cdot 10^{-3}$ & $h=2.5\cdot 10^{-3}$ & $h=10^{-2}$ & $h=5\cdot 10^{-3}$ & $h=2.5\cdot 10^{-3}$\\
                  \hline
$-0.4$ &$1.952e-5$&$9.778e-6$&$4.903e-6$&$2.947e-5$&$1.472e-5$&$7.334e-6$ \\
    \hline
    $-0.375$ &$1.798e-5$&$9.001e-6$&$4.525e-6$&$2.75e-5$&$1.368e-5$&$6.843e-6$ \\
    \hline
    $-0.35$ &$1.68e-5$&$8.424e-6$&$4.199e-6$&$2.526e-5$&$1.264e-5$&$6.359e-6$ \\
    \hline
$-0.325$ &$1.551e-5$&$7.772e-6$&$3.876e-6$&$2.313e-5$&$1.166e-5$&$5.882e-6$ \\
    \hline
$-0.3$ &$1.426e-5$&$7.142e-6$&$3.542e-6$&$2.125e-5$&$1.065e-5$&$5.306e-6$ \\
    \hline
$-0.275$ &$1.317e-5$&$6.527e-6$&$3.25e-6$&$1.943e-5$&$9.76e-6$&$4.906e-6$ \\
    \hline
$-0.25$ &$1.178e-5$&$5.87e-6$&$2.919e-6$&$1.757e-5$&$8.801e-6$&$4.421e-6$ \\
    \hline
$-0.225$ &$1.086e-5$&$5.364e-6$&$2.698e-6$&$1.586e-5$&$7.898e-6$&$3.961e-6$ \\
    \hline
$-0.2$ &$1.023e-5$&$5.116e-6$&$2.573e-6$&$1.57e-5$&$7.808e-6$&$3.902e-6$ \\
    \hline
$-0.175$ &$1.e-5$&$4.994e-6$&$2.519e-6$&$1.559e-5$&$7.714e-6$&$3.83e-6$ \\
    \hline
      \end{tabular}
\end{table}

\begin{table}[!ht]
\centering
\caption{The errors of  a weighted FEM   in the norm of the space ${\bf W}^1_{2, \nu} (\Omega) $
(rotation form), $\omega=9 \pi/8.$
 }
\medskip
  \begin{tabular}{|c|c|c|c|c|c|c|}
    \hline
                & \multicolumn{3}{c|}{ $\delta=0.0123, \nu=1.8$ } & \multicolumn{3}{c|}{ $\delta=0.0127, \nu=1.7$}  \\
       \hline
$\nu^{\ast}$        & $h=10^{-2}$  & $h=5\cdot 10^{-3}$ & $h=2.5\cdot 10^{-3}$ & $h=10^{-2}$ & $h=5\cdot 10^{-3}$ & $h=2.5\cdot 10^{-3}$\\
                  \hline
$-0.3$ &$1.54e-5$&$7.658e-6$&$3.825e-6$&$2.616e-5$&$1.296e-5$&$6.498e-6$ \\
    \hline
    $-0.275$ &$1.418e-5$&$7.071e-6$&$3.519e-6$&$2.376e-5$&$1.192e-5$&$5.937e-6$ \\
    \hline
    $-0.25$ &$1.288e-5$&$6.499e-6$&$3.206e-6$&$2.147e-5$&$1.069e-5$&$5.354e-6$ \\
    \hline
$-0.225$ &$1.166e-5$&$5.847e-6$&$2.896e-6$&$1.919e-5$&$9.501e-6$&$4.763e-6$ \\
    \hline
$-0.2$ &$1.053e-5$&$5.263e-6$&$2.621e-6$&$1.715e-5$&$8.534e-6$&$4.256e-6$ \\
    \hline
$-0.175$ &$9.831e-6$&$4.927e-6$&$2.472e-6$&$1.544e-5$&$7.685e-6$&$3.86e-6$ \\
    \hline
$-0.15$ &$8.542e-6$&$4.247e-6$&$2.127e-6$&$1.347e-5$&$6.716e-6$&$3.344e-6$ \\
    \hline
$-0.125$ &$8.03e-6$&$3.932e-6$&$1.954e-6$&$1.244e-5$&$6.253e-6$&$3.137e-6$ \\
    \hline
      \end{tabular}
\end{table}

\begin{table}[!t]
	\centering\small
	\caption{ The proportion of nodes of the standard FEM, the error in which in absolute value does not exceed
given ones $ \Xi_n $ for a nonlinear problem
(7)-(9) (convective form).}
	\vspace*{2mm}
	\begin{tabular}{|c|c|c|c|c|}
		\hline
		 $\omega_m$ & $\Xi_n$  &  $h=10^{-2}$ &$h=5\cdot 10^{-3}$  &$h=2.5 \cdot 10^{-3}$\\\hline
 $3\pi/2$ & $1\cdot10^{-6}$  & $0.072$ & $0.114$ & $0.178$\\
& $2.5\cdot10^{-6}$ & $0.123$ & $0.195$ & $0.308$\\\hline
 $5\pi/4$ & $5\cdot10^{-7}$  & $0.1$ & $0.322$ &$0.472$\\
& $1\cdot10^{-6}$ & $0.2$ & $0.42$  & $0.607$\\\hline
 $9\pi/8$ & $2.5\cdot10^{-7}$ & $0.116$ & $0.346$ & $0.57$\\
& $5\cdot10^{-7}$ & $0.215$ & $0.576$ & $0.832$\\\hline
			\end{tabular}
\label{tabl_8}
\end{table}

\begin{table}[!t]
	\centering\small
	\caption{The proportion of nodes of the weighted FEM, the error in which in absolute value does not exceed
given ones $ \Xi_n$ for a nonlinear problem
(7)-(9) (convective form). }
	\vspace*{2mm}
	\begin{tabular}{|c|c|c|c|c|c|c|c|}
		\hline
		 $\omega_m$ & $\nu$ &  $\delta$ &  $\nu^{\ast}$ & $\Xi_n$  & $h=10^{-2}$ &$h=5\cdot 10^{-3}$  &$h=2.5 \cdot 10^{-3}$\\\hline
 $3\pi/2$ & 2.0& $0.0127$ & $-0.275$ & $1\cdot10^{-6}$  & $0.104$ & $0.248$ & $0.461$\\
& & & & $2.5\cdot10^{-6}$ & $0.179$ & $0.426$ &$0.696$\\\hline
 $5\pi/4$ & 2.0& $0.0124$ & $-0.15$  & $5\cdot10^{-7}$  & $0.176$ & $0.399$ & $0.761$\\
& & & & $1\cdot10^{-6}$ & $0.276$ & $0.598$ &$0.903$\\\hline
 $9\pi/8$ & 2.0& $0.0127$ & $-0.025$ & $2.5\cdot10^{-7}$ & $0.208$  &$0.435$  & $0.775$\\
& & & & $5\cdot10^{-7}$ & $0.291$ & $0.694$  & $0.937$ \\\hline
			\end{tabular}
\label{tabl_9}
\end{table}

\paragraph{Conclusions.}

A numerical method has been developed for solving stationary nonlinear Navier-Stokes
 equations in convective and rotation
  form
in the domain
with an reentrant corner greater than $ \pi$. The results of the numerical experiments
of the test problems showed
advantages of the method over the classical approach.
The approximate solution for the classical FEM converges to
exact solution of the problem
(7)-(9) (velocity vector) with order $ 0.54 $ for $ \omega_1=1.5 \pi, $
  with order $ 0.67 $ for $ \omega_2=1.25 \pi $ and
  with order $ 0.8 $ for $ \omega_3 =1.125 \pi $ with respect to the grid step $ h $
in the norm of the space $ {\bf W}^1_2(\Omega_m). $
The approximate solution for the weighted FEM (velocity vector) converges to
exact solution of the problem   (7)-(9) for both forms with velocity $ {\cal O}(h) $
 in the norm of the space
$ {\bf W}^ 1_ {2, \nu} (\Omega_m) $ for all values $ \omega_m, m = 1,2,3, $
of the reentrant corner $ \omega$   (see Tables 1-7).

  It is found that free parameters of the method $ \delta $ and $ \nu $
  do not depend on the value of the entering corner $ \omega, $ $ \delta $ is the radius of the neighborhood
  the singularity point is proportional to the grid step $ h$, and the exponent $ \nu $ of the weight
function in the definition of $ R_{\nu}$ generalized solution --- 2.
Parameters $ \nu^{\ast} $ and $ \mu^{\ast} $
(exponents of weight functions in a finite element basis)
  take equal negative values.

The approach is simple to implement and does not require refining the mesh
 in the vicinity of the reentrant corner.
 The method allows us with good accuracy
  find a solution near a singular point and does not allow,
unlike classical approaches, errors spread
to the inner part of the computational domain, where the solution is sufficiently smooth
(see Tables 8-9).

\paragraph{Acknowledgments.}

The reported study was supported by Russian Science Foundation, project No. 21-11-00039,
https://rscf.ru/en/project/21-11-00039/. The results were obtained using the equipment of SRC ”Far Eastern
Computing Resource” IACP FEB RAS (https://cc.dvo.ru).

\end{document}